\documentclass[12pt]{article}

\usepackage[top=1in, bottom=1in, left=1.25in, right=1.25in]{geometry}
\usepackage{amsfonts, amssymb, graphics, amsmath, amsthm, upgreek, comment, tikz, scalefnt, graphicx, caption, float, subcaption, booktabs, cite, amscd,authblk}
\usepackage[romanian,english]{babel}
\usepackage[mathscr]{euscript}
\usepackage[normalem]{ulem}

\numberwithin{equation}{subsection}

\newtheorem{thm}{Theorem}

\newtheorem{lem}{Lemma}[thm]

\theoremstyle{definition}
\newtheorem*{rmk}{Remark}
\newtheorem*{defin}{Definition}
\newtheorem*{no}{Notation}
\newtheorem*{rk}{Remark}

\numberwithin{lem}{section}

\DeclareMathOperator{\R}{\mathbb{R}}
\DeclareMathOperator{\Z}{\mathbb{Z}}

\DeclareMathOperator{\dist}{dist}

\def\v{\vspace{0.1in}}
\def\vv{\vspace{0.2in}}

\def\vecu{\mathbf{u}}
\def\uu{\mathbf{u}}
\def\zz{\mathbf{z}}
\def\zhat{\mathbf{\hat{z}}}
\def\vecx{\mathbf{x}}
\def\x{\mathbf{x}}
\def\vecf{\mathbf{f}}
\def\xun{(\vecx_n,\vecu_n)}
\def\xu0{(\vecx_0,\vecu_0)}
\def\xn{\mathbf{x}_n}
\def\un{\uu_n}

\def\bs{\bar{s}}
\def\bsn{\bar{s}_n}
\def\bso{\bar{s}_0}

\def\C{\mathcal{C}}
\def\G{\mathcal{G}}
\def\S{\mathcal{S}}
\def\H{\mathcal{H}}
\def\P{\mathcal{P}}

\usepackage{anysize}

\begin{document}

\title{Knots Connected by Wide Ribbons}
\author[1]{Susan C. Brooks}
\affil[1]{Department of Mathematics and Philosophy, Western Illinois University - Quad Cities}
\author[2]{Oguz Durumeric}
\affil[2]{Department of Mathematics, University of Iowa}
\author{Jonathan Simon}
\affil[3]{Department of Mathematics, University of Iowa}
\date{\today}
\maketitle


\begin{abstract}
A ribbon is, intuitively, a smooth mapping of an annulus $S^1 \times I$ in 3-space having constant width $\varepsilon$.  This can be formalized as a triple $(x,\varepsilon, \mathbf{u})$ where $x$ is smooth curve in 3-space and $\mathbf{u}$ is a unit vector field based along $x$.  In the 1960s and 1970s, G. C\u{a}lug\u{a}reanu, G. H. White, and F. B. Fuller proved relationships between the geometry and topology of thin ribbons, in particular the ``Link = Twist + Writhe" theorem that has been applied to help understand properties of double-stranded DNA.  Although ribbons of small width have been studied extensively, it appears that less is known about ribbons of large width whose images (even via a smooth map) can be singular or self-intersecting..  

Suppose $K$ is a smoothly embedded knot in $\mathbb{R}^3$.  Given a regular parameterization $\mathbf{x}(s)$, and a smooth unit vector field $\mathbf{u}(s)$ based along $K$, we may define a ribbon of width $R$ associated to $\mathbf{x}$ and $\mathbf{u}$ as the set of all points $\mathbf{x}(s) + r\mathbf{u}(s)$,  $r \in [0,R]$.  For large $R$, these  wide ribbons typically have self-intersections.  In this paper, we analyze how the knot type of the outer ribbon edge $\mathbf{x}(s) + R\mathbf{u}(s)$ relates to that of the original knot $K$.   

We show that, generically, there is an eventual limiting knot type of the outer ribbon edge as $R$ gets arbitrary large. We prove that this eventual knot type is one of only finitely many possibilities which depend just on the vector field $\mathbf{u}$.    However, the particular knot type within the finite set depends on  the parameterized curves $\mathbf{x}(s)$, $\mathbf{u}(s)$, and their interactions.  Finally, we show how to control the curves and their parameterizations so that given two knot types $K_1$ and $K_2$, we can find a smooth ribbon of constant width connecting curves of these two knot types.\\
 
\noindent AMSC: 57M25, 53A04, 53A05
\end{abstract}


\section{Introduction}
\label{sec:intro}

A closed ``ribbon" is a smooth mapping (or the image set) of an annulus, $S^1 \times [0,1]$ into $\R^3$, where the sets ${s} \times [0,1]$ are mapped to line segments all of the same length. To avoid degenerate situations, we assume the mapping is an embedding on $S^1 \times \{0\}$ and  think of a ribbon as being a smooth closed curve with the surface growing out of it following a vector field emanating from the curve. Note we are {\em not} assuming the line segments emanating from $S^1 \times \{0\}$ are orthogonal to the curve, but the restriction to  {\em constant width} is essential in our study

A thin ribbon does not self-intersect, and the ribbon itself gives an isotopy (which extends to an ambient isotopy)  of the two boundary curves. The geometry of thin ribbons has proven to be important in the study of double-stranded DNA  (see e.g. \cite{DS92, ER96}). The key is the   ``link = twist + writhe"  theorem (\cite{GC61,FF71,JW69}).  

\v
We are led to several questions about wide ribbons, the first being how the knot types of the boundary curves can be related. Wide ribbons generally do intersect themselves, and the boundary curves can be of different knot types.  As a ribbon is allowed to grow arbitrarily wide, does the knot type of the outer boundary curve stabilize to something we can predict? Does the restriction to constant width limit which knot types can be connected to which others? 


\subsection{Examples: Knot type may, or may not, change.}
\begin{figure}[h]
\center{
\includegraphics[width=1.5in]{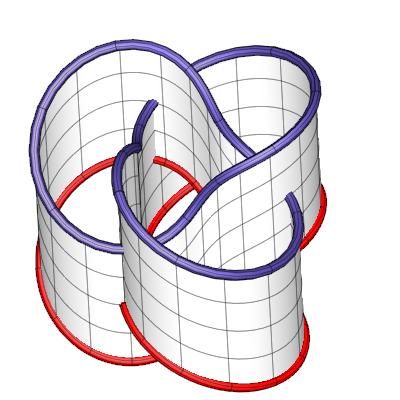}
\caption{Wide ribbon with identical boundary knots}
\label{fig:ParallelKnots}
}
\end{figure}

Having the ribbon self-intersect does not force the outer boundary curve to cross itself. For example, start with any smooth knot and let the vector field be just one constant vector. The ribbon will eventually intersect itself, but the outer boundary curve remains a rigid copy of the original knot, as in  Figure \ref{fig:ParallelKnots}.
 
\begin{figure}[h]
\center{
$\vcenter{\hbox{\includegraphics[width=1.2in] {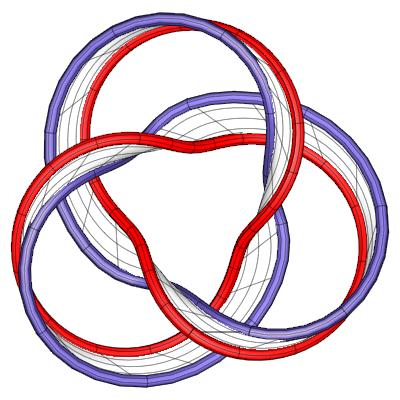} }}$
$\vcenter{\hbox{\includegraphics[width=1.3in] {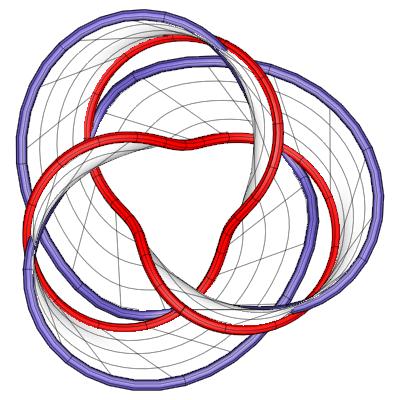} }}$ 
$\vcenter{\hbox{\includegraphics[width=1.55in] {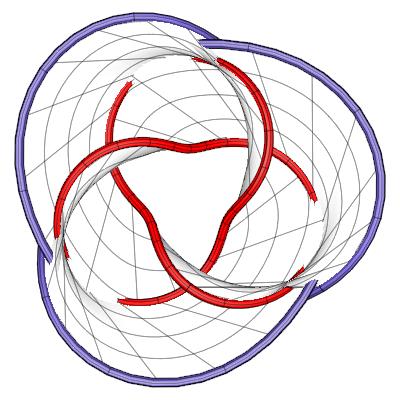} }}$
$\vcenter{\hbox{\includegraphics[width=1.75in] {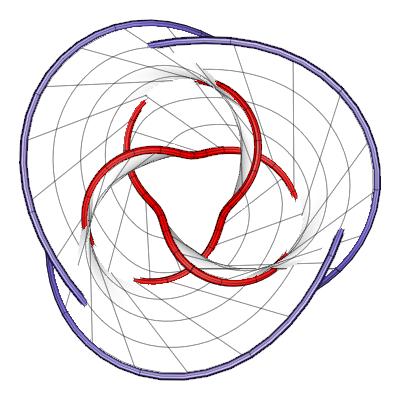} }}$
\caption{Right-hand trefoil flips to left-hand as ribbon gets wider}
\label{fig:changetypes}
}
\end{figure}

On the other hand, with more more general vector fields,  the outer boundary knot can change, as in Figure \ref{fig:changetypes}.


\subsection{Main Results}

It is conceivable that the outer ribbon edge is self-intersecting for all sufficiently large widths, or for some unbounded sequence. However, we show that {\em generically}, this does not happen: in general, as a ribbon grows beyond some width, the outer ribbon edge does not cross itself any more, so the outer ribbon edge eventually stabilizes to a fixed embedded knot type (Theorem \ref{thm:nogoalposts}).  

In Section \ref{sec:goalpost}, we introduce a geometric condition ({\em``no goal posts"})  that ensures this eventual stabilization.  In Section \ref{sec:stability}, we show the condition of having no goal posts is a generic property of ribbons in an appropriate topology (Theorem \ref{generic}).  In Section \ref{sec:bound}, we show that for such non-degenerate ribbons, the limiting knot type is one of finitely many knot types that can be listed from the starting data:  in particular, this set of possible knot types is determined by the vector field along which the ribbon is expanding and not the starting boundary curve (Theorem \ref{thm:2kknottype}).  

The particular choice from the finite set depends on interplay between the initial curve, the vector field, and how they are parameterized.  In Section \ref{sec:constructrib}, we control the curves and parameterizations to show that given any two knot types, we can construct a ribbon of constant width whose boundary curves represent those knot types (Theorem \ref{thm:constructribbon}). 


\subsection{Definitions and assumptions}\label{GeneralAssumptions}
We assume all maps are $\mathcal{C}^{1}$-smooth.  Let $D$ denote $\R/\Z \cong S^1$.
When we talk about a smooth, closed curve, we mean, in particular, that the curve is smoothly closed.  Let $X$ be a smoothly embedded closed curve in $\mathbb{R}^3$, with regular parameterization  $\mathbf{x}: D \rightarrow \mathbb{R}^3$.  Without loss of generality, assume that $\mathbf{x}$ is of unit speed.  

Let $\uu: D \rightarrow S^2$ be a regular smooth closed curve on the unit sphere $S^2$. 
We make the following additional assumption, analogous to the usual knot theory definition of ``regular projection": {\em $\uu$, as a curve on $S^2$, has only transversal self-intersections, and there are no triple points.}  Since the domain $D$ is compact, there can be only finitely many pairs $s \neq \bar{s}$ where $\vecu(s) = \vecu(\bar{s})$.  If we view the vector $\vecu(s)$ as being based at the point $\vecx(s)$, we can think of the function $\vecu$ as a smooth unit vector field along the curve $X$.

\begin{defin}  The \emph{ribbon of width $R$} associated to $\mathbf{x}$ and $\mathbf{u}$ is defined to be $Y= \{\mathbf{x}(s) + r\mathbf{u}(s) \;|\; s \in D \mbox{ and } r \in [0, R]\}$.  The \emph{outer ribbon edge} $Y_{R}$  of $Y $ is the set of points $\mathbf{y}_{R}(s) = \mathbf{x}(s) + R\mathbf{u}(s)$.
\end{defin}

\begin{rmk}
Our notion of {\em width} is more general than some others, because we do not assume  $\mathbf{u}$ is perpendicular to $\mathbf{x}$. We allow $\mathbf{u}(s) \cdot \mathbf{x}'(s)$ to vary. 
\end{rmk}




Intuitively, one might expect that the knot type of the outer ribbon edge $Y_R$ should stabilize for sufficiently large $R$.  The property we define in the next section is a potential obstruction to such stabilization. In Section \ref{sec:goalpost}, we show that the absence of this property, hence the desired stabilization, is in fact generic. So almost all ribbons have eventually constant knot type of the outer curve.



\section{The Goal Post Property}
\label{sec:goalpost}

\begin{defin}
Points $\mathbf{x}(s)$ and $\mathbf{x}(\bar{s})$ along the knot are said to have the \emph{goal post property} with respect to the vector field $\mathbf{u}$ if 
\begin{itemize}
\item $\mathbf{x}(s) \neq \mathbf{x}(\bar{s})$
\item $\mathbf{u}(s) = \mathbf{u}(\bar{s})$ 
\item $\mathbf{u}(s) \cdot \left(\mathbf{x}(s)-\mathbf{x}(\bar{s})\right) = 0$.
\end{itemize}
\end{defin}

The pair $(\mathbf{x},\mathbf{u})$ has the goal post property if there exists such a pair of points. We will show that if the outer ribbon edge $Y_R$ crosses itself for arbitrarily large $R$ then $(\mathbf{x},\mathbf{u})$ has the goal post property.

\vv
Suppose  $s$ and  $\bar{s}$ are distinct parameter values for which the outer ribbon edge intersects itself at some positive width $R$; that is, $\mathbf{y}_{R}(s) = \mathbf{y}_{R}(\bar{s})$.  Then
\begin{equation}\label{EqnForCrossing}
\mathbf{x}(s) - \mathbf{x}(\bar{s}) = R(\mathbf{u}(\bar{s}) - \mathbf{u}(s)). 
\end{equation}

If this happens for a given pair $(s, \bar{s})$ and two widths, $R_1, R_2$, then we have
$$ R_1(\mathbf{u}(\bar{s}) - \mathbf{u}(s)) = \mathbf{x}(s) - \mathbf{x}(\bar{s}) = R_2(\mathbf{u}(\bar{s}) - \mathbf{u}(s))\;.$$
So $\vecu(s)=\vecu(\bar{s})$ or $R_1=R_2$. But we cannot have $\vecu(s)=\vecu(\bar{s})$ since (equation \ref{EqnForCrossing}) this would imply $\vecx(s)=\vecx(\bar{s})$. We thus have a well-defined function defined on those pairs $(s, \bar{s})$ where $\vecu$ crosses itself.

\begin{no}
For each pair of distinct parameters $(s, \bar{s})$, either the rays emanating from $\vecx(s)$ and $\vecx(\bar{s})$ never meet, or there is a single width, which we denote $R(s,\bar{s})$, at which they cross. 
\end{no}

We also will make use of the following lemma, which is obtained by applying the Mean Value Theorem in each coordinate.

\begin{lem}
\label{lem:seqderiv}
If $\vecf:\R \rightarrow \mathbb{R}^k$ is $\mathcal{C}^1$ and we have two sequences $\{s_n\}$ and $\{\bar{s}_n\}$  \em{($s_n \neq \bar{s}_n$)} converging to the same limit,  $s_n \rightarrow s_0$ and $\bar{s}_n \rightarrow s_0$, then
$$\lim_{n \rightarrow \infty} \frac{\vecf(s_n) - \vecf(\bar{s}_n)}{s_n - \bar{s}_n} = \vecf'(s_0).$$
\end{lem}

We can now establish our first theorem: {\em If there are no goal posts, then the outer ribbon edge eventually stabilizes.}

\begin{thm}
\label{thm:nogoalposts}
Suppose we have a smooth closed curve $X$, with parameterization $\vecx(s)$ and unit vector field $\vecu(s)$ satisfying the conditions specified in Section (\ref{GeneralAssumptions}). 
Let $\mathcal{R}$ denote  the set of all widths at which the outer curve $Y_R$ fails to be embedded; 
that is, $\mathcal{R} = \{R(s,\bar{s}) \; |\; s \neq \bar{s} \mbox{ and } 
\mathbf{x}(s) - \mathbf{x}(\bar{s}) =
 R(s,\bar{s})(\mathbf{u}(\bar{s}) - \mathbf{u}(s))\}$. 
 
 If there are no goal posts, then the set $\mathcal{R}$ is bounded and the knots $Y_R$ are isotopic to each other for all $R>\mathop{sup} \mathcal{R}$.
\end{thm}

\begin{proof} If $\mathcal{R}=\emptyset$ then all curves $Y_R$ are isotopic to $X$.
Suppose $\mathcal{R}$ is nonempty and unbounded. Then we can find convergent sequences $s_n \to s_0$ and $\bar{s}_n \to \bar{s}_0$ ($s_n \neq \bar{s}_n$) such that $\lim R(s_n, \bar{s}_n) = \infty$. We will show  $s_0 = \bar{s}_0$ implies $\vecu'(s_0)=\mathbf{0}$ (which contradicts our regularity condition on $\vecu$), and $s_0 \neq \bar{s}_0$ implies the existence of goal posts. Let $R_n$ denote $R(s_n, \bar{s}_n)$.

Suppose first that $s_0 = \bar{s}_0$.
From equation (\ref{EqnForCrossing}), we have
\begin{equation}
\label{eq:seqbound}
\frac{||\mathbf{u}(s_n) - \mathbf{u}(\bar{s}_n)||}{||\mathbf{x}(s_n) - \mathbf{x}(\bar{s}_n)||} = \frac{1}{R_n} \longrightarrow 0.
\end{equation}
On the other hand, applying Lemma \ref{lem:seqderiv} separately to $\vecu$ and $\vecx$ , we have
\begin{equation}
\frac{||\mathbf{u}(s_n) - \mathbf{u}(\bar{s}_n)||}{||\mathbf{x}(s_n) - \mathbf{x}(\bar{s}_n)||} \longrightarrow \frac{||\vecu'(s_0)||  }{||\vecx'(s_0)||} \end{equation}
so $||\vecu'(s_0)||=0$.

\v
Now, suppose $s_0 \neq \bar{s}_0$.  From Equation (\ref{eq:seqbound}), we see that
$$\frac{||\mathbf{u}(s_0) - \mathbf{u}(\bar{s}_0)||}{||\mathbf{x}(s_0) - \mathbf{x}(\bar{s}_0)||} = 0\;,$$
so $\mathbf{u}(s_0) = \mathbf{u}(\bar{s}_0)$.

To see that $\mathbf{x}(s_0)$ and $\mathbf{x}(\bar{s}_0)$ have the goal post property, we need only show that $\mathbf{u}(s_0) \cdot (\mathbf{x}(s_0) - \mathbf{x}(\bar{s}_0)) = 0$.  Consider the isosceles triangles whose vertices are $\mathbf{x}(s_n)$, $\mathbf{x}(\bar{s}_n)$ and $\mathbf{x}(s_n) + R_n\mathbf{u}(s_n) = \mathbf{x}(\bar{s}_n) + R_n \mathbf{u}(\bar{s}_n)$ as shown in Figure \ref{fig:isotris}. Since the length of the base edge $|\mathbf{x}(s_n) - \mathbf{x}(\bar{s}_n)|$ is bounded (by the diameter of the knot) and the sides (of length $R_n$) get arbitrarily large, the  base angles $\beta_n$ converge to $\pi/2$. 

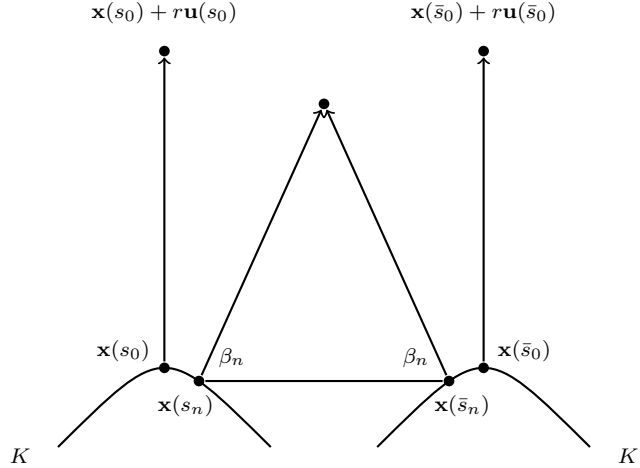
\begin{figure}[H]
\begin{center}
\tikzset{font=\scriptsize}


{\scalefont{0.8}
\begin{tikzpicture}[inner sep=0pt,thick,
        dot/.style={fill=black,circle,minimum size=4pt},scale=3.5]

\draw[thick] (0.2,0) .. controls (0.6,0.4) .. (1.0,0);
\draw[thick] (-0.2,0) .. controls (-0.6,0.4) .. (-1.0,0);

\node[dot](a) at (-0.6,0.3) {};
\node[above left](aa) at (-0.65,0.32){$\mathbf{x}(s_0)$};
\node[dot](b) at (0.6,0.3) {};
\node[above right](bb) at (0.65,0.32){$\mathbf{x}(\bar{s}_0)$};
\node[dot](c) at (-0.47,0.25) {};
\node[below](cc) at (-0.52,0.20){$\mathbf{x}(s_n)$};

\node[dot](d) at (0.47,0.25) {};
\node[below](cc) at (0.52,0.20){$\mathbf{x}(\bar{s}_n)$};
\node[above right](aaa) at (-0.4,0.3){$\beta_n$};
\node[above left](bbb) at (0.4,0.3){$\beta_n$};

\node[dot](e) at (0,1.3) {};

\node[below left](ff) at (-1.1,0){$K$};
\node[below right](gg) at (1.1,0){$K$};

\draw (c) -- (d);
\draw[->] (c) -- (e);
\draw[->] (d) -- (e);

\node[dot](h) at (0.6,1.5) {};
\node[dot](i) at (-0.6,1.5) {};

\draw[->] (b) -- (h);
\draw[->] (a) -- (i);

\node[above](ii) at (-0.6,1.6) {$\mathbf{x}(s_0) + r\mathbf{u}(s_0)$};
\node[above](hh) at (0.6,1.6) {$\mathbf{x}(\bar{s}_0) + r\mathbf{u}(\bar{s}_0)$};

\end{tikzpicture}}


\end{center}
\caption{Isosceles triangle formed near a goal post.}
\label{fig:isotris}
\end{figure}

Now, if we assume that no distinct pair of points $\mathbf{x}(s)$ and $\mathbf{x}(\bar{s})$ has the goal post property, then by the above argument, there exists an $R^* \in \mathbb{R}$ such that $R^* = \mathop{sup} \mathcal{R}< \infty.$  In other words, there does not exist a distinct pair of parameter values $s$ and $\bar{s}$ such that $R(s,\bar{s}) > R^*$.  Hence, $\mathbf{y}_R$ serves as an ambient isotopy of $\mathbf{y}_{R_1}$ to $\mathbf{y}_{R_2}$ for $R^* < R_1 \leq R \leq R_2 < \infty$, implying that the knot type of $Y_{R}$ is unique for $R > R^*$.
\end{proof}


\section{Generic Stabilization of the Outer Ribbon Edge}
\label{sec:stability}

\vv
In this section, we show that stabilization of the outer ribbon edge is a generic property of ribbons, in the sense that within the space of all pairings $(\mathbf{x},\uu)$, the subspace of those having no goal posts is open and dense. Since the term ``ribbon" includes a specified width $R$, we will refer to a pair $(\vecx, \vecu)$ as a {\em ribbon frame}.

\v
Let $\H_1 = \{\vecx:D \to \R^3\}$ and $\H_2 = \{\vecu:D \to S^2\}$ where $\vecx$ and $\vecu$ are $\C^1$ smooth maps.  In each factor, use a $\mathcal{C}^1$ metric:
\begin{eqnarray*}
d_1(\vecx_1, \vecx_2) = \max_{s \in D}(||\vecx_1(s)-\vecx_2(s)||) \; + \; \max_{s \in D}(||\vecx_1'(s)-\vecx_2'(s)||)
\\
 d_2(\vecu_1, \vecu_2) = \max_{s \in D}(||\vecu_1(s)-\vecu_2(s)||) \;+\;\max_{s \in D}(||\vecu_1'(s)-\vecu_2'(s)||)
 \end{eqnarray*}

Let $\H = \H_1 \times \H_2$ and use $d_1+d_2$ as the metric on $\H$: 
$$d((\vecx_1,\vecu_1), (\vecx_2,\vecu_2)) = d_1(\vecx_1,\vecx_2)+d_2(\vecu_1,\vecu_2)\;.$$

Let  $\G$ be the subset of $\H$ consisting of all  $(\vecx, \vecu)$ satisfying the various conditions in Section (\ref{GeneralAssumptions}), and in addition, having no goal posts.  Specifically, the pairs $(\vecx, \vecu) \in \G$ are those where

\begin{enumerate}
\item $\vecx$ and $\vecu$ are regular maps, i.e. $|\vecx'|$ and $|\vecu'|$ are never 0, so both are immersions; \label{xuregular}
\item $\vecx$ is an embedding; \label{xembedded}
\item $\vecu$ has no triple points; \label{notriplepoints}
\item at each double point of $\vecu$, the self-intersection is transversal; \label{uselftransversal}
\item (no goal posts) whenever   $\vecu(s_1) = \vecu(s_2)$,  ($s_1 \neq s_2$), we have $\vecu(s_1) \cdot (\vecx(s_1)-\vecx(s_2)) \neq 0$. \label{nogoalposts}
\end{enumerate}

We wish to show that $\G$ is open and dense.  To that end, we begin by proving two general lemmas regarding sequences of functions.

\begin{lem} \label{lem:fnsn}
Suppose $P,Q$ are metric spaces, $P$ compact, and $f_n:P \to Q$ where $(f_n)$ is a sequence of continuous functions converging uniformly to a continuous map $f:P \to Q$. If $s_n$ is a sequence of points in $P$ converging to $s_0 \in P$, then $\lim_nf_n(s_n)=f(s_0)$.
\end{lem}

\begin{proof}[Proof of Lemma \ref{lem:fnsn}]
From the triangle inequality, $\dist(f(s_0),f_n(s_n)) \leq \dist(f(s_0),f_n(s_0)) + dist(f_n(s_0),f_n(s_n))\;.$ The first terms converge to $0$ because $f_n \rightarrow f$ at each point. The uniformly convergent sequence $(f_n)$, with compact domain $P$, is equicontinuous (converse of Arzel\`{a}-Ascoli theorem), and so the second terms also converge to 0. 
\end{proof}

 The next lemma says that if smooth maps $(f_n)$ converge $C^1$ uniformly to an immersion, in particular a locally $1-1$ map, then the functions $f_n$ are eventually locally $1-1$, and in a uniform way.
 
 \begin{rk} We need the derivatives to converge and $f$ an immersion, otherwise $x \to x^3-t^2x$, $t\to 0$,  is an easy counterexample.
 \end{rk}

\begin{lem} \label{lem:LocallyInjective} Suppose $f_n:D=\R/\Z \to \R^3$ where $(f_n)$ is a sequence of smooth maps converging in $C^1$ to an immersion $f:\R/\Z \to \R^3$. Then there exists $\lambda>0$ and index $N$ such that for all $s, \bs \in D$ and all $n \geq N$, if $s \neq \bs$ and $\mathop{dist}(s, \bs) < \lambda$ then $f_n(s) \neq f_n(\bs)$.
\end{lem}

\begin{proof}[Proof of Lemma \ref{lem:LocallyInjective}]

Suppose, to the contrary, that there exists a sequence of parameter pairs $(s_n, \bs_n)$ with $s_n \neq \bs_n$, $\mathop{dist}(s_n, \bs_n) \to 0$, and $f_n(s_n)=f_n(\bs_n)$. Since $D$ is compact, we may assume w.l.o.g. that the parameter sequences converge: $s_n \to s_0$ and $\bs_n \to \bs_0$. Since $\mathop{dist}(s_n, \bs_n) \to 0$, we have that $s_0=\bs_0$.

For each $n$,  since $f_n(s_n)=f_n(\bs_n)$, in particular there is equality in each of the coordinate functions $f_n^1, f_n^2, f_n^3$  of $f_n$. So there exists points $t_n^1, t_n^2, t_n^3$ in the interval between $s_n$ and $\bs_n$ with derivatives
$df_n^1(t_n^1)=0$, $df_n^2(t_n^2)=0$, and $df_n^3(t_n^3)=0$. Since the numbers $t_n^j$ are pinched  between $s_n^j$ and $\bs_n^j$ $(j=1,2,3) $, we know $t_n^j \to s_0^j$.  

Now apply Lemma \ref{lem:fnsn} to each coordinate of the derivatives: 
$t_n^j \to s_0^j$ and $df_n^j \to df^j$ uniformly, so $df_n^j(t_n^j) \to df^j(s_0)$.  Thus we have $\mathbf{0} =\langle df_n^1(t_n^1), df_n^2(t_n^2), df_n^3(t_n^3)\rangle \to \langle df^1(s_0), df^2(s_0), df^3(s_0) \rangle \neq \mathbf{0}$.
\end{proof}


\begin{thm}
The set $\G$ is open and dense in $\H$. \label{generic}
\end{thm}

\begin{proof}
STEP 1: {\bf{$\G$ is open in $\H$}}.

Suppose $(\vecx_n, \vecu_n)$ is a sequence in $\H$ converging to $(\vecx_0, \vecu_0) \in \G$. We want to show that for sufficiently large $n$, $\xun$ eventually satisfies the five properties that characterize $\G$.  The claim is similar in spirit to the stability theorem(s) in \cite{GP74}.

\begin{enumerate}
\item
Since $\x$ and $\uu$ are regular maps, and the domain $D$ is compact, the values of $|\x'|$ and $|\uu'|$ are  bounded away from $0$. With $\C^1$ convergence, the values of $|\xn'|$ and $|\un'|$ are eventually also bounded away from $0$.

\item 
\label{ItemInjective}
We wish to show the maps $\xn$ are $1-1$ for sufficiently large $n$. Suppose, to the contrary, that for infinitely many  $n$, there exist $s_n \neq \bsn$ with $\xn(s_n)=\xn(\bsn)$. Since $D$ is compact, we can extract convergent subsequences and assume $s_n \to s_0$ and $\bsn\to \bso$.

By Lemma \ref{lem:LocallyInjective}, we must have $s_0 \neq \bso$. But then,  Lemma \ref{lem:fnsn} would say $\x(s_0)=\x(\bso)$.

\item 
If infinitely many $\un$ have triple points, then, as in item \ref{ItemInjective}, the map $\vecu$ would have triple points.

\item
\label{transversal}
As before, if infinitely many $\un$ have non-transversal double points, then there are sequences $s_n \to s_0$, $\bsn \to \bso$, with $s_n \neq \bsn$, $\un(s_n)=\un(\bsn)$ and the vectors $\un'(s_n)$ and $\un'(\bsn)$ collinear. Applying Lemma \ref{lem:LocallyInjective} to show $s_0 \neq \bso$, we would have that $\vecu$ has a non-transversal double point.

\item Suppose infinitely many $\xun$ have goalposts. Then argue as in  item \ref{transversal} to conclude that   $(\vecx,\vecu)$ would as well.

\end{enumerate}

STEP 2: {\bf{$\G$ is dense in $\H$}}. 

We are given $(\x, \uu) \in \H$ and want to approximate it with ribbon frames in $\G$. Much of the work is done by classical results (as in \cite{HWhit36}, \cite{JMil07}).  In particular, we can approximate a continuous curve $\x$ in $\R^3$ with smooth embeddings, and we can approximate a continuous curve $\uu$ in $S^2$ with immersions. We can make the first step in our approximation by perturbing $(\x,\uu)$  slightly to assume $(\x,\uu)$ satisfies properties 1 and 2 in the definition of $\G$: Both $\x$ and $\uu$ are regular maps, and $\x$ is an embedding. Note that any sufficiently close first-order approximation of $\uu$ will now preserve being an immersion.
We  further perturb $\uu$ to ensure that $\uu$ has only finitely many pairs of self-intersection (e.g. by stereographic projection to $\R^2$ and then first order Fourier approximation in each coordinate, similar to \cite{AT95, Traut98}); and further adjust $\uu$ to have only double points where the self-intersections are transversal. We are left with the problem of eliminating goal posts.

Suppose $(\mathbf{x}, \mathbf{u})$ has some goal posts.  Let $\S = \{(s_1, \bs_1), \ldots, (s_k, \bs_k)\}$ be the pairs of distinct parameters for which $\uu(s_i) = \uu(\bs_i)$, and let $\S^*$ be those pairs where there are goal posts, i.e. where $\uu(s_i) \cdot (\x(\bs_i)-\x(s_i)) = 0$. Let $F$ be the  map from domain $\tilde{D} = D - \{(s,s)\}$ to $S^2$ given by
$$F(s,t) = \frac{\x(t)-\x(s)}{{\|\x(t)-\x(s)\|}|}\;.$$
Finally, let $\P$ denote the set of all self-crossing points of $\uu$ on $S^2$.

The goal post condition is that a point $\uu(s_i) = \uu(\bs_i)$ lies on the great circle of $S^2$ traced by vectors orthogonal to $F(s_i, \bs_i)$. We have only finitely many such great circles, so almost all uniform rotations of $S^2$ will move the finite set $\P$ off the union of those great circles. Let $\alpha:S^2 \to S^2$ be such a  rotation, as small as we want, and define a new unit vector field along the given knot $X$ by $\tilde{\uu}(s) = \alpha \circ \uu$. Note that $\tilde{\uu}$ has the same set $\S$ of double point parameter pairs as $\uu$, and now, for all these pairs, $\tilde{\uu}(s_i)$ is not orthogonal to $F((s_i, \bs_i))$; i.e. $(\x, \tilde{\uu})$ has no goal posts.

This completes the proof of Theorem \ref{generic}.
\end{proof}


\section{Bounding the Knot Type of the Outer Ribbon Edge}
\label{sec:bound}

For small values of $R$, it is natural to think of the outer ribbon edge as a perturbation of the curve $\mathbf{x}$.  However, as $R$ increases, we gain more insight by viewing the ribbon edge as a perturbation of the spherical curve $\uu$.

\begin{defin}
The \emph{rescaled outer ribbon edge} is  $\zz_t(s) = \frac{1}{R} Y_R = t\mathbf{x}(s) + \uu(s)$ for $t = \frac{1}{R}$.  
{\bf {Since $\zz_t$ is just a scalar multiple of $Y_R$, they are topologically equivalent knots.}} 

To understand the limiting knot type of $Y_R$ as $R \to \infty$, we will analyze $\zz_t$ as $t \to 0$, along with  its normalization (i.e. spherical projection)   
$\zhat_t(s) = \frac{\zz_t(s)}{||\zz_t(s)||}$.  In the following discussions, we often use the phrase ``for $t$ small enough''.  We always require $t>0$, and ``$t$ small enough'' is equivalent to ``$R$ large enough''.

Note that since $||\uu||=1$ and $||\mathbf{x}||$ is bounded, for sufficiently small $t$ we know $||\zz_t|| \neq 0$  and $\zhat_t$ is defined. 
\end{defin}

The functions $\zz_t$ converge uniformly in $\C^1$ to $\uu$ as $t \to 0$. With a little more work, we have the same property for $\zhat_t$. 

\begin{lem}
\label{lem:convergence}
As $t \rightarrow 0$, $\mathbf{\hat{z}_t} \xrightarrow{\C^1} \uu$ uniformly.
\end{lem}

\begin{proof}
By definition,  
$$d_2(\mathbf{\hat{z}_t},\uu) = \max_{s \in D}(||\zhat_t - \uu||) + \max_{s \in D}(||\mathbf{\hat{z}'_t} - \uu'||).$$
Since $\uu$ is unit and $||\mathbf{x}||$ is bounded, $\zz_t $ converges uniformly to $\uu$ and $||\zz_t|| \to 1$.  Similarly, since $||\mathbf{x}'|| = 1$, we have that $\zz'_t $ converges uniformly to $\uu'$.  Now consider $||\mathbf{\hat{z}_t} - \uu||$. For any $s$,  we have the following:
\begin{eqnarray*}
||\zhat_t - \uu|| & \leq & ||\zhat_t - \zz_t|| + ||\zz_t - \uu||\\
& = & \left|\left|\frac{\zz_t}{||\zz_t||} - \zz_t\right|\right| + ||t\mathbf{x}+ \uu - \uu||\\
& = &\left\lvert \;\frac{1}{||\zz_t||} - 1 \;\right\rvert||\zz_t||  + t||\mathbf{x}|| \rightarrow  0
\end{eqnarray*}

\noindent Next, consider $||\zhat'_t - \uu'||$:

$$||\zhat'_t - \uu'|| \leq ||\zz'_t - \zhat'_t|| + ||\zz'_t - \uu'||$$
As noted above, since $\zz'_t $ converges uniformly to $\uu'$,  $||\zz'_t - \uu'|| \rightarrow 0$.  To complete the proof, we need to show  $||\zz'_t - \zhat'_t|| \to 0$.  Using the quotient formula for the derivative of $\zhat_t$ and the fact that $||\zz_t||^2 = \zz_t \cdot \zz_t$ to calculate the derivative of $||\zz_t||$, we have 
\begin{eqnarray*}
||\zz'_t - \zhat'_t|| & = & \left|\left| \zz'_t - \left(\frac{\zz_t}{||\zz_t||}\right)' \right|\right|\\
& = & \left|\left| \zz'_t - \left(\frac{\zz'_t}{||\zz_t||} - \frac{\zz_t \cdot \zz'_t}{||\zz_t||^3}\ \zz_t \right) \right|\right|\\
& = & \left|\left| \frac{\zz'_t}{||\zz_t||}\left(||\zz_t|| -1\right) + \frac{\zz_t \cdot \zz'_t}{||\zz_t||^3}\ \zz_t \right|\right|\\
& \leq & \frac{||\zz'_t||}{||\zz_t||}\;|\left(||\zz_t|| - 1\right)| + \frac{|\zz_t \cdot \zz'_t|}{||\zz_t||^2}
\end{eqnarray*}
Since $||\zz_t|| \rightarrow 1$, $\zz_t' \to \uu'$ which is bounded, and $\zz_t \cdot \zz'_t \rightarrow \uu\cdot \uu' = 0$ since $\uu$ is constant length, both summands converge to 0.
\end{proof}

Next, we show that for $t$ small enough,  the maps $\zhat_t$ look like regular projections of spatial knots into the sphere.

\begin{lem}
There exists $t_0>0$ such that for $t<t_0$: 
\begin{itemize}
\item $\zhat_t$ is an immersion, 
\item each self-intersection of $\zhat_t$ is a transversal double point.
\end{itemize}
\label{lem:ZhatTnice}
\end{lem}

\begin{proof}
From Lemma \ref{lem:convergence},  $||\zhat_t'||$ converges uniformly to $||\uu'||>0$, so the maps $\zhat_t$ are eventually immersions.

For small $t$, the map $D \times [0,t_0) \to S^2$ given by $(s,t) \to \zhat_t(s)$  is a smooth homotopy in $S^2$ between maps $\zhat_t$ and $\uu$. The Transversality Theorem  \cite{GP74} implies that the property that $\uu$ has transversal self-intersections in $S^2$ is a stable property, so  $\zhat_t$ eventually has only transversal self-intersections. 

It remains to show that (eventually) the self-intersections of $\zhat_t$ are only double points.
Suppose there exists a sequence $t_n \rightarrow 0$ such that $\hat{\zz}_{t_n}$ has triple points.  Then there exist distinct parameter values $a_n$, $b_n$, and $c_n$ such that 
$$\hat{\zz}_{t_n}(a_n) = \hat{\zz}_{t_n}(b_n) = \hat{\zz}_{t_n}(c_n).$$
By compactness, there exists convergent subsequences
$$a_n \rightarrow a_0 \hspace{1 cm} b_n \rightarrow b_0 \hspace{1 cm} c_n \rightarrow c_0\;.$$
We have two cases: Either all three are distinct, or some two are equal, say $a_0=b_0$.

Applying Lemma \ref{lem:LocallyInjective} to $\zhat_t \to \uu$, we know that for $t$ sufficiently small, there is a positive lower bound $\lambda$, uniform in $t$, on the distances $|a_n-b_n|$. But if $a_0=b_0$, these distances would have to become arbitrarily small.

If all three $\{a_0, b_0, c_0\}$ are different, then $\uu$ has a triple point, contradicting our initial assumption that $\uu$ has only double points.
\end{proof}

\noindent We can paraphrase the combination of the previous sections, Lemma \ref{lem:convergence}, and Lemma \ref{lem:ZhatTnice}  as follows:
\begin{itemize}
\item {\bf Basic limiting properties of $\zz_t$ and $\zhat_t$:}
 Under the generic assumptions of regularity with no ``goal posts", for sufficiently small $t$, the curves $\zz_t$ are oriented, embedded space curves which are smoothly isotopic to one another via the ribbon and converging uniformly to the oriented spherical curve $\uu$.  Furthermore, the spherical projections $\zhat_t$ are oriented regular curves on the sphere, each having only transversal double point self-intersections, and converging 
$\C^1$-uniformly to $\uu$ as $t \to 0$. 
\end{itemize}

\v
We want to characterize the (single) knot type of the  curves $\zz_t$ as being obtained from $\uu$ by resolving the double points of $\uu$ into over- or under-crossings.  We showed above that the curves $\zhat_t$ look like knot projections; now we want to see which knot.

We claim that for sufficiently small t, the double points of $\zz_t$ occur at essentially the same parameter values as for $\uu$, in the same order, with the same orientations of the curves as for $\uu$. The fact that we have $\C^1$-convergence ensures that the handedness of crossings will agree with the corresponding resolution of $\uu$, so the space curves $\zz_t$ have the same extended Gauss code as a particular resolution of $\uu$.
\vv

We establish the desired relationship between double point parameters of $\zhat_t$ and $\uu$ in several steps, sometimes choosing a tolerance $\delta$ on neighborhoods in $D$ of the double point parameter values of $\uu$ and sometimes making $t$ small enough to force $\zhat_t$ to approximate $\uu$ closely enough. 

\begin{figure}[htbp]
\begin{center}
\includegraphics[width=2.0in] {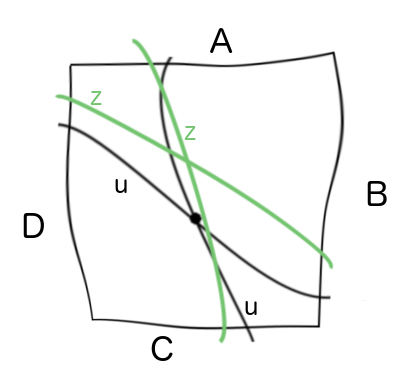}
\end{center}
\caption{Arcs  of $\zhat_t$ are $\C^1$-close to arcs of $\uu$. Spanning across the box neighborhood $\implies$ each double point of $\uu$ has a nearby double point of $\zhat_t$. Crossing at a similar angles $\implies$ double points of $\zhat_t$ are not too close together. }
\label{fig:DoublePointsOfZt}
\end{figure}

Much of the story is told in Figure \ref{fig:DoublePointsOfZt}.  In a box-shaped neighborhood of a double point of $\uu$, we see two arcs of the curve $\uu$ crossing at some angle, with the neighborhood chosen so the arcs span from one ``side" of the box to the opposite side.  (Note that the angle is bounded away from zero due to transversality.) Nearby are two approximating arcs of $\zhat_t$, which cross at a similar angle and follow $\uu$ closely enough to span between the same sides of the box.  The fact that the $\zhat_t$ arcs span across the box implies, by the Jordan Curve Theorem, that they meet somewhere inside the box. At the same time, the fact that their tangent vectors are close to the tangents for $\uu$  prevents $\zhat_t$ from having a second  double point close to the first. These two facts  are the essential ingredients in establishing the desired 1-1 correspondence between double points of $\zhat_t$ and double points of $\uu$.

\subsection{First choice of $\delta$-intervals about the double point parameters of $\uu$}
\label{FirstDelta}
\begin{no} 
Let $\{s_1, \bar{s}_1, s_2, \bar{s}_2, \ldots s_k, \bar{s}_k\}$ be the parameter values for the self-intersection points of $\uu$, where $\uu(s_i) = \uu(\bar{s}_i)$.  Call $s_i, \bar{s}_i$ a {\em matched pair of parameter values}. For $\delta>0$, let $B_\delta(s_i)$ (resp. $\bar{s}_i$) be the interval in the domain $D$ of radius $\delta$ about $s_i$. Call $B_\delta(s_i)$ and $B_\delta(\bar{s}_i)$ a {\em matched pair of intervals}.  Finally, for any $\delta>0$, let $\mathscr{B}$ be the union of the $B_\delta$ neighborhoods of all the double point parameters of $\uu$. 
\end{no}

We first note that by making $t$ small enough, we can ensure that the double point parameters of $\zhat_t$ lie within $\mathscr{B}$. 

\begin{lem}
\label{lem:AllDPPinB}
There exists $t_0>0$ such that for $t<t_0$, all double point parameters of $\zhat_t$ are contained in $\mathscr{B}$.
\end{lem}

\begin{proof}
  
Suppose  there exists a sequence $t_n \rightarrow 0$ where $\hat{\zz}_{t_n}$ has double point parameters $r_n$ and $\bar{r}_n$ such that at least one of the parameter values, call it $r_n$, is not contained in the open set $\mathscr{B}$.  Since $\mathscr{B}$ is open, by choosing convergent subsequences, we may assume $r_n$ converges to $r_0$, which is in the complement of $\mathscr{B}$, and $\bar{r}_n$ converges to $r_\#$ that is somewhere. By Lemma \ref{lem:convergence}, regardless of where $r_\#$ lies, we know $\uu(r_0) = \uu(r_\#)$.

If $r_\# \neq r_0$ then $\uu$ has a double point parameter $r_0$ in the complement of $\mathscr{B}$. If $r_\# = r_0$ then parameters $r_n, \bar{r}_n$ are eventually closer than $\lambda$ from Lemma \ref{lem:LocallyInjective}.
\end{proof}

Since $\uu$ has only finitely many double points, there exists $\delta_1>0$ so that the intervals are disjoint, and in particular:
\begin{itemize}
\item each interval $B_{\delta_1}$ contains {\bf exactly one} double point parameter for $\uu$, and
\item for sufficiently small $t$, all the double point parameters of $\zhat_t$ are contained in the union of these intervals.
\end{itemize}
 Any smaller $\delta$ neighborhoods also isolate the double point parameters of $\uu$; and, at the expense of further decreasing $t$, all of the double point parameters of $\zhat_t$ are contained in the smaller $\mathscr{B}$. So we can choose $\delta_2$ and $t$ small enough that no two matched double point parameters of any $\zhat_t$ are as close as $2\delta_2$. In particular, for $\delta \leq \delta_2$, 
\begin{equation}
\bullet \;\; \text{ no } B_\delta  \text{ can contain a matched pair of double point parameters for any } \zhat_t.
  \label{NoCloseMatchedParameters}
 \end{equation}
 
For future technical reasons, shrink $\delta_2$ if necessary to ensure that 
 \begin{itemize}
 \item the {\bf closed} $B_\delta$ intervals are disjoint.
 \end{itemize}
 
 Again, note that any smaller $t_1$ or $\delta_2$ satisfy the above conditions. This does not yet imply that each interval $B_\delta(s_i)$ (resp. $\bar{s}_i$) contains any  double point parameter for $\zhat_t$, or contains only one - perhaps there are many {\bf{un}}matched parameters inside one $B_\delta$.
 
 
 \subsection{Choose $\delta_3$ and $t_2$ so the $B_\delta$ isolate parameters for $\zhat_t$}
 \label{IntervalsIsolateDoublePoints}
 We know that matched double point parameters of $\zhat_t$ cannot be too close together; we need to establish the same property for {\bf{un}}matched parameters.  This takes a few steps.
 
 We first sharpen the statement that all the double point parameters of $\zhat_t$ are contained in $\mathscr{B}$.
 
\begin{lem} 
 \label{lem:MatchedParameters} 
 For small enough $t$, matched double point parameters of $\zhat_t$ are contained in matched $B_\delta$ intervals. 
\end{lem}

\begin{proof}If not, then there exists $t_n \rightarrow 0$ with matched parameters  $r_n, \bar{r}_n$ such that the $B_\delta$ intervals containing these parameter values are not matched.  As usual, since $D$ is compact, we can extract convergent subsequences: $r_n \rightarrow r_0$ and $\bar{r}_n \rightarrow r_\#$.  Since there are only finitely many $B_\delta$ intervals, we can further extract subsequences so that all $r_n$ are contained in $B_\delta(s_i)$ for one particular $i$. The distance between matched  double point parameters for $\zhat_t$ is bounded away from 0 (recall $\lambda$)  so $r_0 \neq r_\#$.  But then, since $\uu(r_0) = \uu(r_\#)$ by Lemma \ref{lem:convergence}, $r_0$ and $r_\#$ are double point parameters for $\uu$.  Thus, $r_0$ must equal $s_i$ and $r_\# = \bar{s}_i$.  In particular, all $r_n$ are contained in $B_\delta(s_i)$ and all but finitely many $\bar{r}_n$ are contained in the matching interval $B_\delta(\bar{s}_i)$, contradicting the assumption that the pairs $r_n, \bar{r}_n$ are contained in unmatched intervals.
\end{proof}

 We also will use a general property of smooth space curves. 
 
\begin{lem}
Let $\gamma: [a,b] \rightarrow \mathbb{R}^n$ be a $\C^1$ regular curve, let $\mathbf{u}_0$ be a fixed unit vector, and let $\alpha_0$ be some positive angle.  If for each $s \in [a,b]$ we have $\angle(\gamma'(s),\mathbf{u}_0) \leq \alpha_0$, then for each parameter pair $c,d$ with $a \leq c \leq d \leq b$, we have
$$\angle(\gamma(d) - \gamma(c), \mathbf{u}_0) \leq \alpha_0,$$
for $\gamma(c) \neq \gamma(d)$.
\label{lem:chordvector}
\end{lem}

\begin{proof}
Let $\gamma: [a,b] \rightarrow \mathbb{R}^n$ be a $\C^1$ regular curve, and choose parameter values $c$ and $d$ with the property that $a \leq c \leq d \leq b$ and $\gamma(c) \neq \gamma(d)$.  Further, let $s$ denote a parameter value in the domain.  Refer to Figure \ref{fig:chordvector} for a schematic diagram of the curve $\gamma$, an associated chord, and a tangent vector.

\begin{figure}[h]
\begin{center}
\tikzset{every picture/.style={scale=0.7}}
\tikzset{font=\footnotesize}
%
%

\begin{tikzpicture}[inner sep=0pt,thick,
        dot/.style={fill=black,circle,minimum size=4pt}, scale=0.7]

\draw[ultra thick, -] (-4.27,0) -- (3.65,0);
\draw[-] (-4.5,-0.5) .. controls (-3,3) and (-2,2) .. (5,-0.5);
\draw[ultra thick, ->] (-1,1.6) -- (2,0.7);

\node[dot] at (-4.27,0) {};
\node[dot] at (3.65,0) {};
\node[dot] at (-1,1.6) {};

\node at (-5.3, 0) {$\gamma(c)$};
\node at (3.7, 0.8) {$\gamma(d)$};
\node at (0.5, 2) {$\gamma'(s)$};

\end{tikzpicture}

%
\end{center}
\caption{$\gamma$ curve with tangent vector $\gamma'(s)$ and chord $\gamma(d) - \gamma(c)$.}
\label{fig:chordvector}
\end{figure}
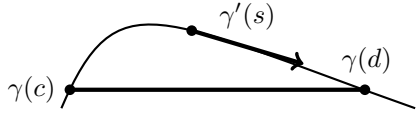

\noindent Fix a unit vector $\mathbf{u}_0$, and let $\alpha_0$ be an acute angle with the property that $\angle(\gamma'(s), \mathbf{u}_0) \leq \alpha_0$ for each $s \in [a,b]$.  Also, let $\theta(s) = \angle(\gamma'(s), \mathbf{u}_0) \leq \alpha_0$.  Since $\mathbf{u}_0$ is unit, we have the following:
\begin{eqnarray*}
\mathbf{u}_0 \cdot (\gamma(d) - \gamma(c)) & = & \mathbf{u}_0 \cdot \int_c^d \gamma'(s)\ ds\\
& = &  \int_c^d  \mathbf{u}_0 \cdot \gamma'(s)\ ds\\
& = & \int_c^d ||\gamma'(s)||\cos\theta(s)\ ds\\
& \geq & \int_c^d ||\gamma'(s)||\cos\alpha_0\ ds\\
& = & \cos \alpha_0 \int_c^d ||\gamma'(s)||\ ds\\
& \geq & \cos\alpha_0\ ||\gamma(d) - \gamma(c)||
\end{eqnarray*}
\noindent But $$\mathbf{u}_0 \cdot (\gamma(d) - \gamma(c)) = ||\gamma(d) - \gamma(c)||\cos\left(\angle(\gamma(d) - \gamma(c), \mathbf{u}_0)\right)$$
Consequently, $$\cos(\angle(\gamma(d) - \gamma(c), \mathbf{u}_0)) \geq \cos\alpha_0 \Longrightarrow \angle(\gamma(d) - \gamma(c), \mathbf{u}_0) \leq \alpha_0$$
\end{proof}

We now can show why unmatched double point parameters of $\zhat_t$ cannot accumulate within one $B_\delta$.

\begin{lem}
There exist   $\delta_3 \leq \delta_2$ and $t_2 \leq t_1$ such that for all  $\delta<\delta_3$, $t<t_2$, no two unmatched double point parameters of $\zhat_t$  are  contained in any one $B_\delta$.
\label{UNmatchedParameters}
\end{lem}

\begin{proof}
With our current values $\delta_2$ and $t_1$, we know that all double point parameters of $\zhat_t$ are contained in the union of the $B_\delta$ intervals and that matched parameters are contained in matched intervals.  Also, if we shrink $\delta$, we can shrink $t$ to preserve these properties.  

At each double point $p_i$ of $\uu$ there are  two arcs of $\uu$ crossing at some positive angle  $\alpha_i$. Let $\alpha_0 = \frac{1}{3}\mathop{min}\{\alpha_1, \ldots, \alpha_k\}$. Since $\zhat_t'$ are uniformly continuous and (Lemma \ref{lem:convergence})  $\zhat_t' \longrightarrow \uu'$ uniformly, there exists $\delta_3>0$ such that for $\delta \leq \delta_3$ and $t$ less than some $t_2$,  if $q$ is a point in $B_\delta(s_i)$ then the angle between $\zhat_t'(q)$ and $\uu'(s_i)$ is less than $\alpha_0$.

We claim that with $\delta < \delta_3$ and $t< t_2$, no $B_\delta(s_i)$ can contain two unmatched double point parameters of $\zhat_t$. Suppose $q,r$ are double point parameters of $\zhat_t$ contained in $B_\delta(s_i)$.  Their matching parameters $\bar{q}, \bar{r}$ are, by Lemma \ref{lem:MatchedParameters}, contained in $B_\delta(\bar{s}_i)$. 
The crossing angle of $\uu$ at $\uu(s_i) = \uu(\bar{s}_i)$ is $\alpha_i$.  By our choice of $\delta_3$ and $t_2$, we have:
\begin{itemize}
\item The angle between $\uu'(s_i)$ and $\zhat_t'(w)$ is less than $\alpha_0$ for each $w \in B_\delta(s_i)$.
\item The angle between $\uu'(\bar{s}_i)$ and $\zhat_t'(w)$ is less than $\alpha_0$ for each $w \in B_\delta(\bar{s}_i)$.
\end{itemize}
Consequently, Lemma \ref{lem:chordvector} ensures the following:
\begin{itemize}
\item The angle between $\uu'(s_i)$ and the chord vector $(\zhat_t(r) - \zhat_t(q))$ is less than $\alpha_0$, and
\item the angle between $\uu'(\bar{s}_i)$ and the chord vector $(\zhat_t(\bar{r}) - \zhat_t(\bar{q}))$ is less than $\alpha_0$
\end{itemize}
But this says that the angle between the two chords is at least $\frac{1}{3}\alpha_i >0$.  On the other hand, if both pairs of parameters are matched, then the chords are identical.
\end{proof}

Once we know that no two double point parameters of any $\zhat_t$ (whether matched or unmatched) lie in a single $B_\delta$, we can say the following:

\begin{itemize}
\item With $\delta<\delta_3$ and $t < t_2$ as above, each interval $B_\delta$ contains {\bf at most one} double point parameter value for a given $\zhat_t$.  (Note we are not yet claiming that $\zhat_t$ {\em has} any double points, much less that the parameter values are close to those for $\uu$; simply that no $B_\delta$ contains two of them for a given $\zhat_t$.)
 \end{itemize}
 
 Again, note that for any smaller choices of $\delta_3$  we can shrink $t_2$ to still satisfy all bulleted properties listed so far.


\subsection{Choose $\delta$ intervals to constrain how  $\uu$ crosses itself and force $\zhat_t$ to have nearby double points}
\label{FindDoublePointsOfZ}
This is another step involving  several steps of local analysis and ``epsilonics". The arguments are similar to the previous section, so we will summarize here.

For each double point $p_i$ of $\uu$, we find a choice of $\delta^i< \delta_3$ so that $\uu$ is especially well-behaved in $B_i = B_{\delta^i}(s_i) \cup B_{\delta^i}(\bar{s}_i)$; we choose $t$ small enough to have $\zhat_t$ follow $\uu$ closely enough to produce double point parameters for $\zhat_t$ in $B_i$.  Then choose $\delta_4$ to be the {\bf maximum}  of these separate $\delta$, and choose $t_3$ the minimum of the respective $t$. This will yield double points parameters of $\zhat_t$ in each $B_{\delta_4}(s_i)$ [resp. $\bar{s}_i$] for all $t<t_3$.

For each double point $p_i = \uu(s_i) = \uu(\bar{s}_i)$ of $\uu$, by referring to the tangent plane $T_{p_i}S^2$, we can find a neighborhood $N(p_i)$ in $S^2$ with the following properties (refer to Figure \ref{fig:DoublePointsOfZt}):
\begin{itemize}
\item $N(p_i)$ in $S^2$  is diffeomorphic to a rectangle with consecutive sides $A, B, C, D$.
\item The various neighborhoods $N(p_i)$ are pairwise disjoint.
\item There exists $0< \delta_4^i \leq \delta_3$ such that the images $\uu(B_{\delta_4^i}(s_i))$ and $\uu(B_{\delta_4^i}(\bar{s}_i))$ are arcs spanning $N(p_i)$ such that one arc connects the interior of side $A$ to the interior of side $C$ and the other arc runs similarly between $B$ and $D$.
\item Each arc $\uu(B_{\delta_4^i}(s_i))$ and $\uu(B_{\delta_4^i}(\bar{s}_i))$ meets the boundary edges of $N(p_i)$ transversally.
\item The derivative $\uu'$ is nearly constant on each  interval $B_{\delta_4^i}(s_i)$ and $B_{\delta_4^i}(\bar{s}_i)$.
\end{itemize} 

For each $p_i$, we obtain the neighborhood $N(p_i)$ by radial projection of an appropriate neighborhood of $p_i$ in the tangent plane $T_{p_i}(S^2)$.

\begin{figure}[h]
\begin{center}
\includegraphics[width=2.0in]{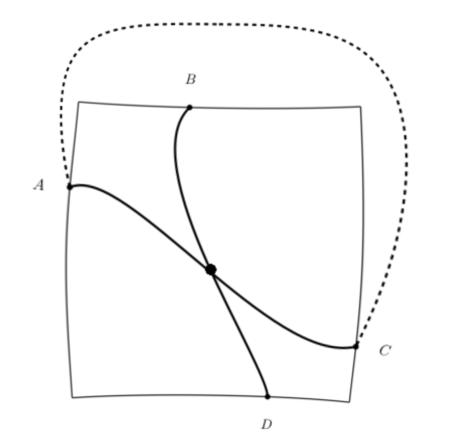}
\end{center}
\caption{Box neighborhood of a double point of $\uu$.}
\label{fig:BoxNbd}
\end{figure}

For each $i$, we can choose $t_3^i$ so that when $t<t_3^i$, the maps $\zhat_t$ are close enough to $\uu$ that  the images $\zhat_t(B_{\delta_4^i}(s_i))$ and $\zhat_t(B_{\delta_4^i}(\bar{s}_i))$ {\bf contain} arcs (the arcs may extend further) in $N(p_i)$ that connect $A$ to $C$ and $B$ to $D$. The Jordan Curve theorem then implies that those arcs must intersect (see Figure \ref{fig:BoxNbd}).
Letting $\delta_4 = \max \{\delta_4^i \} \leq \delta_3$ and  $t_3 = \min \{t_3^i\} \leq t_2$, we have
\begin{itemize}
\item For $t<t_3$, $\zhat_t$ has {\bf at least one} double point parameter in each  $B_{\delta_4}$ interval about a double point parameter of $\uu$.
\item From Section \ref{IntervalsIsolateDoublePoints}, we then have that for $t<t_3$, $\zhat_t$ has {\bf exactly one} double point parameter in each   interval $B_{\delta_3}(s_i)$ [resp $\bar{s}_i$]. And the corresponding double point of $\zhat_t$ lies in the neighborhood $N(p_i)$.
\end{itemize}

We now have intervals about the double point parameters of $\uu$  that isolate the double point parameters of $\uu$, isolate any double point parameters of $\zhat_t$ they happen to contain, AND that each contain exactly one double point parameter of $\zhat_t$, with matched double point parameters of $\zhat_t$ contained in matched $B_\delta$ neighborhoods of the double point parameters of $\uu$.  Thus, the double point parameter values of $\uu$ are in one-to-one correspondence with the double point parameter values of $\zhat_t$.


\subsection{Conclusion: The limiting knot type}
 
We now show that the limiting knot type of $\zz_t$ is one of finitely many choices, which are determined by $\uu$.  Recall that $k$ denotes the number self-intersections of $\uu$.

\begin{thm}
Under our generic assumptions of regularity with {\em no goal posts}, for all sufficiently small $t$, the knot type of $\zz_t$ is constant and is the same as one of the $2^k$ resolutions of $\uu$.
\label{thm:2kknottype}
\end{thm}

\begin{proof}
 For appropriate $\delta$ and sufficiently small $t$, the self-intersections of $\hat{\zz}_t$ are in one-to-one correspondence with the self-intersections of $\uu$ in the strong sense that  matched double point parameters of $\zhat_t$ occur in the  matched $B_{\delta}$ neighborhoods of double point parameters of $\uu$.  Suppose $s_t, \bar{s}_t$ are matched double point parameters of $\zhat_t$. Relative to projection into $S^2$, one of the points $\zhat_t(s_t), \zhat_t(\bar{s}_t)$ lies over the other. Resolve the corresponding crossing of $\uu$ in the same way to obtain an embedded knot $\tilde{\uu}$ with the same  Gauss code as $\zhat_t$. 
 Because $\zhat' \approx \uu'$, when we resolve the crossings of $\uu$ to get the embedded oriented knot $\tilde{\uu}$, the handedness of each crossing of $\tilde{\uu}$ is the same as the handedness of the corresponding crossing of $\zhat_t$. So   $\tilde{\uu}$ and $\zhat_t$ have the same extended Gauss code.
  \end{proof}
  

\section{Constructing Ribbons Between Any Two Knots}
\label{sec:constructrib}

If one begins with a particular knot with parameterization $\mathbf{x}$ and unit vector field $\mathbf{u}$ satisfying our generic assumptions of regularity and no goal posts, Theorem \ref{thm:2kknottype} shows that the outer ribbon edge eventually stabilizes to a resolution of $\mathbf{u}$.  We now show that if we are given two knot types, then it is possible to construct a ribbon frame so that $\mathbf{x}$ is one of the given knot types and the limiting resolution of $\mathbf{u}$ is the other given knot type.


\begin{no}
Recall that $\{s_1, \bar{s}_1, s_2, \bar{s}_2, \ldots, s_k, \bar{s}_k\}$ are the parameter values for the self-intersection points of $\uu$, where $\uu(s_i) = \uu(\bar{s}_i)$.  For $t$ small enough, since the self-intersections of $\uu$ are in one-to-one correspondence with $\zhat_t$, let $\{s_{t,1}, \bar{s}_{t,1}, s_{t,2}, \bar{s}_{t,2}, \ldots, s_{t,k}, \bar{s}_{t,k}\}$ denote the parameter values for the self-intersection points of $\zhat_t$, where $\zhat_t(s_{t,i}) = \zhat_t(\bar{s}_{t,i})$.
\end{no}

\begin{thm}
Given knot types $K_1$ and $K_2$, there exists a ribbon frame $(\mathbf{x},\mathbf{u})$ satisfying the conditions in Section \ref{sec:stability}, where $\mathbf{x}(t)$ defines a knot of type is $K_1$, and the limiting knot type of the outer ribbon edge is type $K_2$.
\label{thm:constructribbon}
\end{thm}

\begin{proof}
We begin by considering a special case, which is illustrative of the general process.\\

\noindent \textbf{Special Case}: $K_1$ is the unknot, and $K_2$ is any knot.\\

By Theorem 3.6 of \cite{DE04}, there exists a smoothly embedded knot $\tilde{\mathbf{u}}$ of type $K_2$ in $\mathbb{R}^3$ with a regular projection into the plane such that there is an arc in the projection which traverses all of the crossings once before traversing any of them a second time.  (Note that such a projection need not be one of minimal crossing number.)  Let $P$ denote the projection mapping into the plane, and let $A$ denote the Hamiltonian arc in $P(\tilde{\mathbf{u}})$.

Fix an orientation and starting point for $A$, and let $\{p_1, p_2, \ldots, p_k\}$ denote the double points of the projection $P(\tilde{\mathbf{u}})$ in the order that they lie along $A$.  For each $j \in \{1,2, \ldots, k\}$, in addition to the $p_j$ label, we will assign a sign, denoted by $\omega(p_j)$, to indicate whether the arc is crossing over or under along $A$.  Let $\omega(p_j) = +1$ denote that $p_j$ is an over-crossing double point, and let $\omega(p_j) = -1$ denote that $p_j$ is an under-crossing double point.  Thus, the double points of $P(\tilde{\mathbf{u}})$ have associated pairs 
$$(p_1, \omega(p_1)), (p_2, \omega(p_2)), \ldots (p_k, \omega(p_k)).$$
See Figure \ref{fig:labelfigeight} for an example of a projection of the figure eight knot as well as a figure of an appropriate labelling of the same projection where the arc $A$ runs between the two circular nodes with the indicated orientation.

\begin{figure}[h]
\begin{center}
   \begin{subfigure}[b]{0.48\textwidth}
                \centering
		\includegraphics[scale=.25]{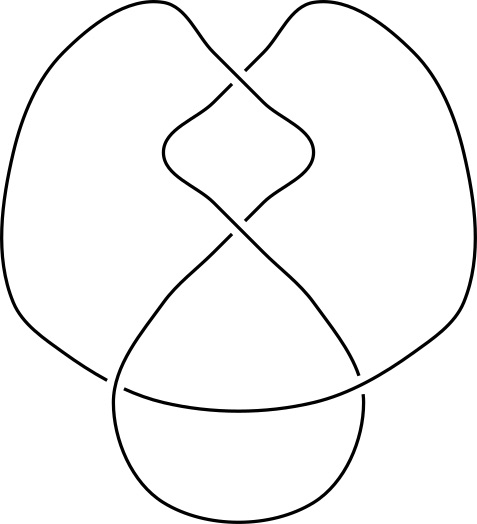}
    \end{subfigure} 
   \begin{subfigure}[b]{0.48\textwidth}
                \centering
		\includegraphics[scale=.25]{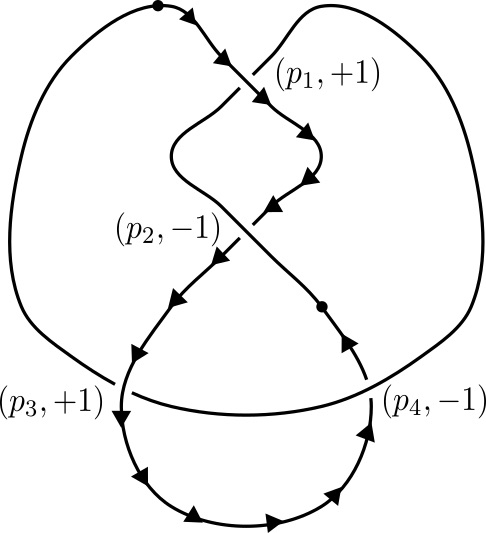}
    \end{subfigure}\end{center}
\caption{A figure eight and its associated labelled diagram.}
\label{fig:labelfigeight}
\end{figure}

Once we traverse through the double point labeled $p_k$ on the knot, we begin to traverse through each crossing a second time, but possibly in a different order.  We will denote these double points (which lie outside of $A$) by the pairs 
$$\left(\bar{p}_{\uptau(1)}, \omega(\bar{p}_{\uptau(1)})\right), \left(\bar{p}_{\uptau(2)}, \omega(\bar{p}_{\uptau(2)})\right), \ldots, \left(\bar{p}_{\uptau(k)}, \omega(\bar{p}_{\uptau(k)})\right),$$
where $(\uptau(1), \uptau(2), \ldots, \uptau(k))$ is a permutation of $(1,2,\ldots,k)$.  We will address these points when we specify parameterizations for the curves that we define to be $\mathbf{x}$ and $\mathbf{u}$.

We wish to separate the over-crossing double points from the under-crossing double points.  To do so, we first perform an ambient isotopy of the plane so that the set of double points $\{p_1, p_2, \ldots, p_k\}$ along $A$ are collinear, i.e. straighten the arc $A$.  We then separate the points by isotoping the over-crossing double points to one side of the line and the under-crossing double points to the other side of the line.  When separating the points, we do so in such a way that we do not introduce any new crossings.  (Note that this is possible since there are only finitely many double points, and an arbitrarily small perturbation is enough to isotope any given point off of the line.). See Figure \ref{fig:linefigeight} for an example of the figure eight knot in Figure \ref{fig:labelfigeight} with collinear double points before and after such an isotopy.

\begin{figure}[h]
\begin{center}
   \begin{subfigure}[b]{0.48\textwidth}
                \centering
		\includegraphics[scale=.35]{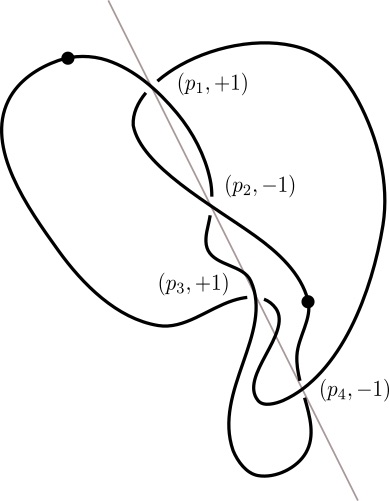}
    \end{subfigure} 
   \begin{subfigure}[b]{0.48\textwidth}
                \centering
		\includegraphics[scale=.35]{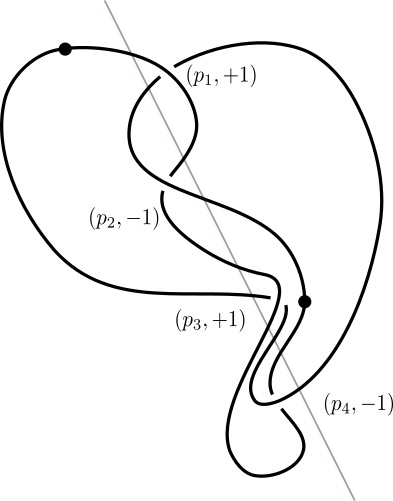}
    \end{subfigure}
\end{center}
\caption{A figure eight with collinear double points before and after isotopy.}
\label{fig:linefigeight}
\end{figure}

Once the set of double points $\{p_1, p_2, \ldots, p_k\}$ has been separated according to their sign within the plane, we enclose each of the two groups within disks.  Then we smoothly isotope the plane to $S^2 \setminus \{(0, -1, 0)\}$ so that the disks map to small polar caps (e.g. with polar angle less than $\frac{\pi}{12}$) at the north and south pole of $S^2$.  We choose the isotopy so that the over-crossing points lie within the disk at the north pole while the under-crossing points lie within the disk at the south pole.  The straight line  originally containing the double points is mapped to the equator.  We define this isotopic version of $P(\tilde{\mathbf{u}})$ on $S^2$ to be $\mathbf{u}$.

Now, we construct an appropriate curve to represent $\mathbf{x}$.   After defining the geometric curve $\mathbf{x}$, we will adjust the parameterization for $\mathbf{x}$ appropriately to create the association $\mathbf{x}(s) \leftrightarrow \mathbf{u}(s)$.  Since $K_1$ is the unknot, let $\mathbf{x}: D \rightarrow S^2$ be the smooth arclength parameterization of the great circle where $\mathbf{x}(0) = \mathbf{x}(2\pi) = (0,0,1)$, $\mathbf{x}\left(\frac{\pi}{2}\right) = (1,0,0)$, and $\mathbf{x}(\pi) = (0,0,-1)$.  We now reparameterize $\mathbf{x}$ to control where $\mathbf{x}$ maps the double point parameters of $\mathbf{u}$:  
$$0 < s_1 < s_2 < \cdots < s_k < \bar{s}_{\uptau(1)} < \bar{s}_{\uptau(2)} < \cdots < \bar{s}_{\uptau(k)} < 2\pi$$
	\begin{enumerate}
	\item Compress $\mathbf{x}\ | [0,s_k]$ so that the set $\mathbf{x}([0,s_k])$ is contained in the small cap at the north pole.
	\item Stretch $\mathbf{x}\ | [s_k, \bar{s}_{\uptau(1)}]$ so that $\{\mathbf{x}(\bar{s}_{\uptau(1)})\}$ is contained in the small cap at the south pole.
	\item Compress $\mathbf{x}\ | [\bar{s}_{\uptau(1)},\bar{s}_{\uptau(k)}]$ so that the set $\{\mathbf{x}([\bar{s}_{\uptau(1)},\bar{s}_{\uptau(k)}])\}$ is contained in the small cap at the south pole.
	\item Stretch $\mathbf{x}\ | [\bar{s}_{\uptau(k)}, 2\pi]$ to complete the great circle.
	\end{enumerate}

To emphasize the relation between the planar and spherical projections, we use the label $p_j$ to denote $\mathbf{u}(s_j)$ on $S^2$, and likewise let $\bar{p}_{\uptau(j)}$ denote $\mathbf{u}(\bar{s}_{\uptau(j)})$.  Our choice of parameterizations for $\mathbf{x}$ and $\mathbf{u}$ allows us to control the over-crossing and under-crossing pattern of our outer ribbon edge $\mathbf{z}_t$ in order to achieve the desired knot type.  Indeed, suppose that $\mathbf{u}(s_j) = \mathbf{u}(\bar{s}_j)$ is near the north pole.  This means that the double point $p_j$ along $A$ was on an over-crossing strand.  Also recall that $\mathbf{x}(s_j)$ is near the north pole while $\mathbf{x}(\bar{s}_j)$ is near the south pole.  Our choice of parameterization for $\mathbf{x}$, together with the fact that $\zz_t$ is uniformly close to $\uu$ implies that $\uu(s_{t,j}) \cdot \mathbf{x}(s_{t,j}) \approx 1$ while $\uu(\bar{s}_{t,j}) \cdot \mathbf{x}(\bar{s}_{t,j}) \approx -1$. So $||\mathbf{z}_t(s_{t,j})|| = ||\mathbf{u}(s_{t,j}) + t\mathbf{x}(s_{t,j})||  > 1$  while $ ||\mathbf{z}_t(\bar{s}_{t,j})|| = ||\mathbf{u}(\bar{s}_{t,j}) + t\mathbf{x}(\bar{s}_{t,j})||  < 1$.  That is,
$$||\mathbf{z}_t(s_{t,j})|| > ||\mathbf{z}_t(\bar{s}_{t,j})||$$
Conversely, if $\mathbf{u}(s_j) = \mathbf{u}(\bar{s}_j)$ is near the south pole, then $||\mathbf{z}_t(s_{t,j})|| < 1$ and $||\mathbf{z}_t(\bar{s}_{t,j})|| > 1$ so that
$$||\mathbf{z}_t(s_{t,j})|| < ||\mathbf{z}_t(\bar{s}_{t,j})||$$
Thus, $\mathbf{z}_t$ has the same over-crossing and under-crossing configuration as $\tilde{\mathbf{u}}$ and, therefore, stabilizes to $K_2$.  (It is also important to note that our choice of parameterizations and placements of double points ensures that $\mathbf{u}(s_i) \cdot \mathbf{x}(s_i) \approx - \mathbf{u}(\bar{s}_i) \cdot \mathbf{x}(\bar{s}_i) \approx \pm 1$.  This guarantees that the ribbon frame $(\mathbf{x}, \mathbf{u})$ does not have the goal post property.)\\


\noindent \textbf{General Case}: \emph{$K_1$ and $K_2$ represent arbitrary knot types.}\\

Let $\mathbf{x}$ be a smoothly embedded curve in $\mathbb{R}^3$ whose knot type is that of $K_1$.  We can proceed as in the special case above by controlling the behavior of $\mathbf{x}$ near the self-intersections of $\mathbf{u}$, which lie near the polar caps.  As such, we can isotope the curve so that $\mathbf{x}$ follows the great circle through $(0,0,\pm 1)$ and $(-1,0,0)$ except for a sufficiently small ball centered at $(1,0,0)$ containing the ``interesting'' part of the knot.  For an example of a suitable $\mathbf{x}$ curve whose knot type is the trefoil, see Figure \ref{fig:partref}.\\

\begin{figure}[h]
\begin{center}
\tikzset{every picture/.style={scale=0.8}}
\tikzset{font=\footnotesize}



\begin{tikzpicture}[y=0.80pt,x=0.80pt,yscale=-1, inner sep=0pt, outer sep=0pt, scale=.15,dot/.style={fill=black,circle,minimum size=4pt}]
\begin{scope}
  \path[draw=black,line join=miter,line cap=butt,miter limit=4.00,line
    width=1.200pt] (375.0000,302.3622) .. controls (375.0000,302.3622) and
    (370.2923,343.8908) .. (368.2589,393.9584)(367.6294,458.0886) .. controls
    (368.3484,483.8265) and (370.4979,508.2911) .. (375.0000,527.3622) .. controls
    (387.8436,581.7685) and (437.1564,622.9559) .. (450.0000,677.3622) .. controls
    (455.7438,701.6934) and (467.6777,734.6845) .. (450.0000,752.3622) .. controls
    (437.3884,764.9738) and (419.1865,773.0867) ..
    (399.3826,776.7010)(335.3347,772.8958) .. controls (321.5887,768.4874) and
    (309.2807,761.6429) .. (300.0000,752.3622) .. controls (229.2893,681.6515) and
    (229.2893,523.0729) .. (300.0000,452.3622) .. controls (335.3553,417.0068) and
    (414.6447,417.0068) .. (450.0000,452.3622) .. controls (467.6777,470.0398) and
    (455.7438,503.0310) .. (450.0000,527.3622) .. controls (446.2415,543.2835) and
    (439.3599,558.0728) .. (431.1830,572.3926)(393.8170,632.3318) .. controls
    (385.6401,646.6516) and (378.7585,661.4409) .. (375.0000,677.3622) .. controls
    (357.7685,750.3558) and (375.0000,902.3622) .. (375.0000,902.3622);

\draw (200,550) circle (700pt);

\node at (500,302) {$\mathbf{x}$};
\node at (-200,-400) {$S^2$};


\draw[thin]  (-470,-10) .. controls (-23,150) and (423,150) .. (870,-10);
\draw[thin, dashed]  (-470,-10) .. controls (-23,-210) and (423,-210) .. (870,-10);

\draw[thin]  (-470,1110) .. controls (-23,1270) and (423,1270) .. (870,1110);
\draw[thin, dashed]  (-470,1110) .. controls (-23,910) and (423,910) .. (870,1110);

\draw[very thick] (375,302.3622) ..  controls (360, -300) and (200, -325) .. (200,-320);
\draw[very thick] (375,902.3622) .. controls (360, 1400) and (200, 1425) .. (200, 1420);

\draw[very thick, dashed] (200,-325) .. controls (0, 258) and (0, 841) .. (200, 1425);

\end{scope}

\end{tikzpicture}

\end{center}
\caption{Appropriate $\mathbf{x}$ curve near $S^2$ representing the knot type of a trefoil.}
\label{fig:partref}
\end{figure}
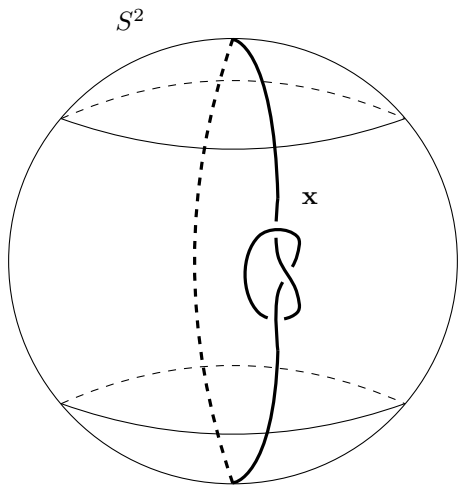

Once the curve defined by $\mathbf{x}$ has been fixed, we may proceed by defining the parameter values $\{s_1, s_2, \ldots, s_k\}$ and $\{\bar{s}_{\uptau(1)}, \bar{s}_{\uptau(2)}, \ldots, \bar{s}_{\uptau(k)}\}$ as explained in the previous case, which depend on $\mathbf{u}$ alone and not $\mathbf{x}$.  The remainder of the arguments in the previous case apply to such a curve.
\end{proof}

\newpage

\bibliography{bibliography}{}

\begin{thebibliography}{10}

\bibitem{GC61}
G.~C\u{a}lug\u{a}reanu.
\newblock Sur les classes d{'}isotopie des n\oe uds tridimensionnels et leurs
  invariants.
\newblock {\em Czech. Math. J.}, 11:588--625, 1961.

\bibitem{DE04}
Y.~Diao, C.~Ernst, and X.~Yu.
\newblock Hamiltonian knot projections and lengths of thick knots.
\newblock {\em Topology and its Applications}, 136:7--36, 2004.

\bibitem{FF71}
F.~Brock Fuller.
\newblock The writhing number of a space curve.
\newblock {\em Proc. Nat. Acad. Sci. USA}, 68(4):815--819, 1971.

\bibitem{GP74}
Victor Guillemin and Alan Pollack.
\newblock {\em Differential Topology}.
\newblock Prentice-Hall, Inc., Englewood Cliffs, NJ, 1974.

\bibitem{JMil07}
J.~Milnor.
\newblock {\em Collected Papers of John Milnor: III Differential Topology}.
\newblock American Mathematical Society, Providence, RI, 2007.

\bibitem{DS92}
De~Witt~L. Sumners.
\newblock Knot theory and {D}{N}{A}.
\newblock In DeWitt~L. Sumners, editor, {\em New Scientific Applications of
  Geometry and Topology, Proceedings of Symposia in Applied Mathematics},
  volume~45, pages 39--72. Am. Math. Soc., 1992.

\bibitem{AT95}
Aaron Trautwein.
\newblock {\em Harmonic Knots}.
\newblock PhD thesis, Department of Mathematics, University of Iowa, Iowa City,
  Iowa, 1995.

\bibitem{Traut98}
Aaron Trautwein.
\newblock An introduction to harmonic knots.
\newblock In Andrzej Stasiak, Vsevolod Katritch, and Louis~Hirsch Kauffman,
  editors, {\em Ideal Knots}, volume~19, pages 353--373. World Scientific
  Publishing Co., 1998.
\newblock Series on Knots and Everything.

\bibitem{ER96}
E.~J.~Janse van Rensburg, Enzo Orlandini, De~Witt Sumners, M.~Carla Tesi, and
  Stuart~G. Whittington.
\newblock Topology and geometry of biopolymers.
\newblock In Jill~P. Mesirov, Klaus Schulten, and De~Witt Sumners, editors,
  {\em Mathematical Approaches to Biomolecular Structure and Dynamics},
  volume~82, pages 21--38. Springer-Verlag Publ., 1996.
\newblock IMA Volumes in Mathematics and its Applications.

\bibitem{JW69}
J.~H. White.
\newblock Self-linking and the {G}auss integral in higher dimensions.
\newblock {\em Am. J. Math.}, 91:693--728, 1969.

\bibitem{HWhit36}
H.~Whitney.
\newblock Differentiable manifolds.
\newblock {\em Annals of Mathematics}, 37:645--680, 1936.

\end{thebibliography}
\bibliographystyle{plain}

\end{document}